\documentclass[letterpaper, 10 pt, conference]{ieeeconf}
\IEEEoverridecommandlockouts

\usepackage{amsmath}
\usepackage{graphicx}
\usepackage{cite}
\usepackage{authblk}
\usepackage{amssymb}
\usepackage{subfigure}
\usepackage{multirow}
\usepackage{color}
         
\newtheorem{theorem}{Theorem}

\newtheorem{proposition}{Proposition}

\title{Active Target Defense Differential Game with a Fast Defender}
\author{Eloy Garcia, David W. Casbeer, and Meir Pachter  
\thanks{A preliminary version of this manuscript has been submitted to the 2015 American Control Conference.}
\thanks{E. Garcia is a contractor (Infoscitex Corp.) with the Control Science Center of Excellence, Air Force Research Laboratory, Wright-Patterson AFB, OH 45433. \ttfamily{elgarcia@infoscitex.com}}
\thanks{D. Casbeer is with the Control Science Center of Excellence, Air Force Research Laboratory, Wright-Patterson AFB, OH 45433. \ttfamily{david.casbeer@us.af.mil}}
\thanks{M. Pachter is with the Department of Electrical Engineering, Air Force Institute of Technology, Wright-Patterson AFB, OH 45433. \ttfamily{meir.pachter@afit.edu}}
}

\begin{document}
\maketitle 
\begin{abstract}
This paper addresses the active target defense differential game where an Attacker missile pursues a Target aircraft. A Defender missile is fired by the Target's wingman in order to intercept the Attacker before it reaches the aircraft. Thus, a team is formed by the Target and the Defender which cooperate to maximize the distance between the Target aircraft and the point where the Attacker missile is intercepted by the Defender missile, while the Attacker tries to minimize said distance. The results shown here extend previous work. We consider here the case where the Defender is faster than the Attacker. The solution to this differential game provides optimal heading angles for the Target and the Defender team to maximize the terminal separation between Target and Attacker and it also provides the optimal heading angle for the Attacker to minimize the said distance. 
\end{abstract}

\section{Introduction} \label{sec:intro}
Pursuit-evasion scenarios involving multiple agents represent important and challenging types of problems in aerospace, control, and robotics. They are also useful in order to analyze biologically inspired behaviors. For instance, the paper \cite{Scott13} addressed a scenario where two evaders employ coordinated strategies to evade a single pursuer, but also to keep them close to each other. The authors of \cite{Oyler14} discussed a multi-player pursuit-evasion game with line segment obstacles labeled as the Prey, Protector, and Predator Game. 
Dominance regions were provided for each agent in order to solve the game, that is, to determine if the Protector is able to rescue the Prey before the Predator captures it.  
A different approach to address pursuit-evasion games with several pursuers in order to capture an evader within a bounded domain is based on dynamic Voronoi diagrams, as in \cite{Huang11} and \cite{Bakolas10}.   

\addtolength{\parskip}{0mm}
\addtolength{\abovedisplayskip}{-1mm}
\addtolength{\belowdisplayskip}{-1mm}
\addtolength{\dblfloatsep}{-3mm}
\addtolength{\dbltextfloatsep}{-3mm}

An scenario of active target defense including three agents, the Target ($T$), the Defender ($D$), and the Attacker ($A$), has been analyzed in the context of cooperative optimal control \cite{Boyell76}, \cite{Boyell80}. 
Indeed, sensing capabilities of missiles and aircraft allow for implementation of complex pursuit and evasion strategies \cite{Zarchan97}, \cite{Siouris04}, and recent work has proposed different guidance laws for the agents $A$ and $D$. 
In \cite{ratnoo11} the authors addressed the case where the Defender implements Command to the Line of Sight (CLOS) guidance to pursue the Attacker which requires the Defender to have at least the same speed as the Attacker.
A different guidance law for the Target-Attacker-Defender (TAD) scenario was given by Yamasaki \textit{et.al.} \cite{Yamasaki10}, \cite{Yamasaki13}. These authors investigated an interception method called Triangle Guidance (TG), where the objective is to command the defending missile to be on the line-of-sight between the attacking missile and the aircraft for all time while the aircraft follows some predetermined trajectory.
The authors show, through simulations, that TG provides better performance in terms of Defender control effort than a number of variants of Proportional Navigation (PN) guidance laws, that is, when the Defender uses PN to pursue the Attacker instead of TG. These approaches constrain and limit the level of cooperation between the Target and the Defender by implementing Defender guidance laws without regard to the Target's trajectory. 

The papers \cite{rubinsky12}, \cite{rubinsky13} presented an analysis of the end-game TAD scenario based on the Attacker/Target miss distance for a \textit{non-cooperative} Target/Defender. The authors develop linearization-based Attacker maneuvers in order to evade the Defender and continue pursuing the Target. 

Different types of cooperation have been recently proposed in \cite{Perelman11}, \cite{Rusnak05}, \cite{Rusnak11}, \cite{Ratnoo12}, \cite{Shima11}, \cite{Shaferman10}, \cite{Prokopov13} for the TAD scenario.
In these papers the Target represents an aircraft trying to evade a missile homing on it. The Defender represents another missile launched by the aircraft (or a wingman) in order to intercept and destroy the Attacker in order to guarantee the survival of the aircraft. 
Thus, in \cite{Rusnak11} optimal policies (lateral acceleration for each agent including the Attacker) are provided for the case of an aggressive Defender, that is, the Defender has a definite maneuverability advantage. A linear quadratic optimal control problem is posed where the Defender's control effort weight is driven to zero to increase its aggressiveness. Reference \cite{Ratnoo12} provided a game theoretical analysis of the TAD problem using different guidance laws for both the Attacker and the Defender. 
The cooperative strategies in \cite{Shima11} allow for a maneuverability disadvantage for the Defender with respect to the Attacker and the results show that the optimal Target maneuver is either constant or arbitrary. 
In the recent paper \cite{Prokopov13} the authors analyze different types of cooperation assuming the Attacker is oblivious of the Defender and its guidance law is known. Two different one-way cooperation strategies were discussed: when the Defender acts independently, the Target knows its future behavior and cooperates with the Defender, and vice versa. Two-way cooperation where both Target and Defender communicate continuously to exchange their states and controls is also addressed, and it is shown to have a better performance than the other types of cooperation - as expected.

Our preliminary work \cite{Garcia14}, \cite{Garcia15} considered the cases when the Attacker implements typical guidance laws of Pure Pursuit (PP) and PN, respectively. In these papers, the Target-Defender team solves an \textit{optimal control} problem that returns the optimal strategy for the $T-D$ team so that $D$ intercepts the Attacker and at the same time the separation between Target and Attacker at the instant of interception of $A$ by $D$ is maximized. 

In this paper the active target defense scenario is modeled as a zero-sum three-agent pursuit-evasion differential game. The two-agent team consists of a Target and a Defender who cooperate; the Attacker is the opposition. The goal of the Attacker is to capture the Target while the Target tries to evade the Attacker and avoid capture. The Target cooperates with the Defender which pursues and tries to intercept the Attacker before the latter captures the Target. Cooperation between the Target and the Defender is such that the Defender will capture the Attacker before the latter reaches the Target. In this differential game the Attacker also solves an optimal control problem in order to minimize the final separation between itself and the Target. Assuming that the Attacker knows the position of the Defender, this strategy provides better performance for the Attacker than using PP or PN. From the Attacker's point of view, it is better to bring the Defender-Attacker interception point closer to the Target's position (and hopefully produce some damage), even though the Attacker is then captured by the Defender. The present paper extends the results in \cite{Pachter14Allerton} where it was assumed that both missiles, the Attacker and the Defender, have the same speed. Here, we extend the analysis of this differential game to include the operationally relevant case where the Attacker and the Defender missiles have different speeds; the focus of this paper is on the case where the Defender is faster than the Attacker. This scenario is more complex than the previously considered particular case of same speeds. Here, we derive the optimal strategies for each one of the three agents. In addition, given a Defender-Attacker speed ratio we provide the critical Target/Attacker speed ratio to guarantee its survival.  

We also obtain the analytical solutions of the differential game, and give special attention to the case where the Target starts closer to the Attacker than to the Defender. For this scenario we also provide the critical minimal speed of the Target for it to avoid capture; that is, when the Target starts closer to the Attacker than to the Defender, its speed must be bounded from below; otherwise the Target will be captured by the Attacker before the Defender can get in the way of the Attacker and intercept it. 

The paper is organized as follows. Section \ref{sec:Problem} states the active target defense differential game. 
Section \ref{sec:Num} provides a numerical method to solve the differential game. 
The minimum Target speed ratio to evade capture by the Attacker is given in Section \ref{sec:Minspeed}. 
Analytical solutions of the differential game are provided in Section \ref{sec:analysis}. Examples are given in Section \ref{sec:examples} and concluding remarks are made in Section \ref{sec:concl}.

\section{Differential Game} \label{sec:Problem}
\begin{figure}
	\begin{center}
		\includegraphics[width=8.4cm,height=5.5cm,trim=1.7cm 1.6cm 1.7cm .5cm]{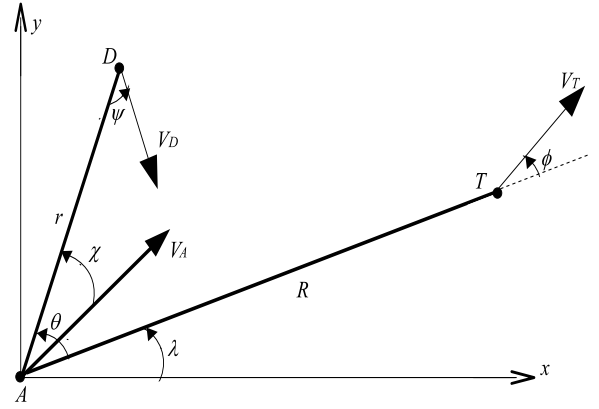}
	\caption{Reduced state space}
	\label{fig:DG}
	\end{center}
\end{figure}

The target defense differential game is illustrated in Fig. \ref{fig:DG}. The speeds of the Target, Attacker, and Defender are denoted by $V_T$, $V_A$, and $V_D$, respectively, which are assumed to be constant. The agents have ``simple motion" a la Isaacs.
The dynamics of the three vehicles in the realistic game space are given by:
\begin{align}
	\dot{x}_T&=V_T\cos\hat{\phi}, \ \ \ \ \ \ \ \ \ \ \ \dot{y}_T=V_T\sin\hat{\phi}     \label{eq:fixedT}  \\   
  \dot{x}_A&=V_A\cos\hat{\chi}, \ \ \ \ \ \ \ \ \ \ \  \dot{y}_A=V_A\sin\hat{\chi}   \label{eq:fixedADG} \\
	\dot{x}_D&=V_D\cos\hat{\psi}, \ \ \ \ \ \ \ \ \ \ \: \dot{y}_D=V_D\sin\hat{\psi}  	\label{eq:fixedD}
\end{align}
where the headings of the Target, the Attacker, and the Defender are, respectively, $\hat{\phi}=\phi+\lambda$, $\hat{\chi}=\lambda\!+\!\theta\!-\!\chi$, and $\hat{\psi}=\psi+\theta+\lambda-\pi$.

The variables $R$ and $r$ represent the separation between the Attacker and the Target and between the Attacker and the Defender, respectively. In this game the Attacker pursues the Target and tries to capture it. The Target and the Defender cooperate
in order for the Defender to intercept the Attacker before the latter captures the Target. Thus, the Target-Defender team search for a cooperative optimal strategy to maximize $R(t_f)$ which represents the distance between the Target and the Attacker at the time instant $t_f$ of the Defender capturing the Attacker. The Attacker will search for its corresponding optimal strategy in order to minimize $R(t_f)$.

Define the speed ratio problem parameter $\alpha=V_T/V_A$. In general, we have that the Attacker missile is faster than the Target aircraft, so $\alpha <1$. 
Let us define the speed ratio $\beta=V_D/V_A$. When the Defender is faster than the Attacker we have that $\beta>1$.

\section{Numerical Solution} \label{sec:Num}
In this section, the corresponding dynamics of the three-agent engagement will be modeled using the reduced state space formed by the ranges $R$ and $r$, and by the angle between them, denoted by $\theta$ see Fig. \ref{fig:DG}. The objective of the Target-Defender team is to determine their optimal heading angles $\phi$ and $\psi$ in this reduced state space such that the distance $R(t_f)$ is maximized at the time instant $t_f$ where the separation $r(t_f)=r_c$, where $r_c$ denotes the Defender's capture radius. The interception time $t_f$ is free. The objective of the Attacker is to determine its optimal heading angle, denoted by $\chi$, such that the distance $R(t_f)$ is minimized. Note that the relative heading angles can be easily transformed to heading angles with respect to the fixed coordinate axis $x$ using the line of sight angle from the Attacker to the Target, denoted by $\lambda$. The use of the reduced state space provides a compact representation of the dynamics of this three-agent differential game.

The (normalized with respect to the speed $V_A$) dynamics in the reduced state space are   
\begin{align}
	\dot{R}&=\alpha \cos\phi - \cos(\theta-\chi), \ \ \ \ \ \ \  R(t_0)=R_0 \label{eq:dynR_DG} \\   
 	\dot{r}&=-\cos\chi - \beta\cos\psi, \ \ \ \ \ \ \ \ \ \ \ \: r(t_0)=r_0  \\
	\dot{\theta}&=-\frac{\alpha}{R}\sin\phi + \frac{1}{R}\sin(\theta-\chi) \nonumber  \\
	 &~~ \ - \frac{\beta}{r}\sin\psi + \frac{1}{r}\sin\chi,  \ \ \ \ \ \ \ \ \ \theta(t_0)=\theta_0 \label{eq:dynD_DG}  
\end{align}
for $0\leq t \leq t_f$.

The objective of the Target-Defender team is to maximize the separation between the Target and the Attacker at the interception time $R(t_f)$, where the terminal time $t_f$ is free, such that $r(t_f)=r_c$. The objective of the Attacker is to minimize the same distance $R(t_f)$. This can be expressed as
\begin{align}
    \max_{\phi,\psi}  \min_\chi	J=\int_{t_0}^{t_f}\dot{R}dt. 	\label{eq:costDG}
\end{align}
Then, the Hamiltonian is given by (where the Target-Defender team aims at minimizing $-J$ and the Attacker aims at maximizing $-J$, for convenience of notation of the solutions):
\begin{align}
\left.
	\begin{array}{l l}
	H\!\!\!\!\!\!&=\cos(\theta-\chi)-\alpha\cos\phi \\
	&~~+(\alpha \cos\phi - \cos(\theta-\chi))\lambda_R \\
	&~~-\left(\cos\chi + \beta\cos\psi\right)\lambda_r  \\      
 	 &~~+ (-\frac{\alpha}{R}\sin\phi \!+\! \frac{1}{R}\sin(\theta\!-\!\chi) \!-\! \frac{\beta}{r}\sin\psi \!+\! \frac{1}{r}\sin\chi)\lambda_\theta  	\label{eq:HamiltonianDG}
	\end{array}  \right.
\end{align}
and the co-state dynamics are given by:
\begin{align}
	\dot{\lambda_R}&=\frac{\lambda_\theta}{R^2}\left(\sin(\theta-\chi) - \alpha\sin\phi\right) \label{eq:costateRDG}  \\ 	
	\dot{\lambda_r}&=\frac{\lambda_\theta}{r^2}\left(\sin\chi - \beta\sin\psi\right)   \\
  \dot{\lambda_\theta}&=(1-\lambda_R)\sin(\theta-\chi) - \frac{\lambda_\theta}{R}\cos(\theta-\chi).  \label{eq:costateDG}
\end{align}
The terminal conditions for this free terminal time problem are as follows. The terminal state $r(t_f)$ is fixed and equal to $r_c$. Because the terminal states $R(t_f)$, and $\theta(t_f)$ are free, we have $\lambda_R(t_f)=\lambda_\theta(t_f)=0$. The final terminal condition for optimality for this problem requires that $H(x^*(t_f),u^*(t_f),\lambda^*(t_f),t_f)=0$. In summary, the terminal conditions are:
\begin{align}
\left.
	\begin{array}{c c}
	  r(t_f)=r_c  \\
	  \lambda_R(t_f)\!\!=0  \;\;	\\
	  \lambda_\theta(t_f)\!=0  \;\;	\\
	  \alpha^2\! + \!2\big(\alpha\beta\! +\! \cos\theta(t_f)\big)\lambda_r(t_f) + (\beta^2-1)\lambda^2_r(t_f) -1\!=0.\;  
	\end{array}   \right. \label{eq:DGtc}
\end{align}

\begin{proposition} \label{th:DG}
The Target and Defender optimal control headings that maximize the separation between the Target and the Attacker and achieve $r(t_f)=r_c$ are given by
\begin{align}
	  \sin\psi^*&=\frac{\lambda_\theta}{r\sqrt{\lambda_r^2 + \lambda_\theta^2/r^2}}   \label{eq:DGspsi}  \\
		\cos\psi^*&=\frac{\lambda_r}{\sqrt{\lambda_r^2 + \lambda_\theta^2/r^2}} 	     \label{eq:DGcpsi}
\end{align}

\begin{align}
	  \sin\phi^*&=\frac{\lambda_\theta}{R\sqrt{(1\!-\!\lambda_R)^2\! +\! \lambda_\theta^2/R^2}}  \label{eq:DGsphi}   \\
		\cos\phi^*&=\frac{1-\lambda_R}{\sqrt{(1-\lambda_R)^2 + \lambda_\theta^2/R^2}}. 	     \label{eq:DGcphi}
\end{align}
The Attacker optimal control heading that minimizes the separation between itself and the Target at $t=t_f$ is given by
\begin{align}
	  \sin\chi^*&=\frac{\chi_s}
		{\sqrt{\chi_s^2+\chi_c^2}}   \label{eq:DGschi}   \\
		\cos\chi^*&=\frac{\chi_c}
		{\sqrt{\chi_s^2+\chi_c^2}}  \label{eq:DGcchi}
\end{align}
where
$\chi_s=(1-\lambda_R)\sin\theta-\frac{\lambda_\theta}{R}\cos\theta+\frac{\lambda_\theta}{r}$ and $\chi_c=(1-\lambda_R)\cos\theta+\frac{\lambda_\theta}{R}\sin\theta-\lambda_r$.
\end{proposition}
\textit{Proof}. 
In order to find the optimal heading angle equations involving $\psi^*$ we solve for this variable by differentiating the Hamiltonian \eqref{eq:HamiltonianDG} in $\psi$ and setting the derivative to $0$
\begin{align}
  \frac{\partial H}{\partial \psi}=\beta\lambda_r \sin\psi - \frac{\beta}{r}\lambda_\theta \cos\psi=0. \label{eq:pHpsiLOS}
\end{align}
Using the trigonometric identity $\cos^2\psi=1-\sin^2 \psi$ we can write \eqref{eq:pHpsiLOS} as
\begin{align}
  \sin^2\psi= \frac{\lambda_\theta^2}{r^2 (\lambda_r^2 + \lambda_\theta^2/r^2)} \nonumber
\end{align}
and \eqref{eq:DGspsi} follows. The expression \eqref{eq:DGcpsi} is found in a similar way by letting $\sin^2\psi=1-\cos^2 \psi$ in \eqref{eq:pHpsiLOS}.
We can compute the second partial derivative of the Hamiltonian with respect to $\psi$ to show that this solution minimizes the cost $-J$. Doing so we obtain  
\begin{align}
\left.
	\begin{array}{l l}
  \frac{\partial^2 H}{\partial \psi^2}&=\beta\lambda_r \cos\psi + \frac{\beta}{r}\lambda_\theta \sin\psi \\
	&=\frac{\beta\lambda_r^2}{\sqrt{\lambda_r^2 + \lambda_\theta^2/r^2}} + \frac{\beta\lambda_\theta^2}{r^2\sqrt{\lambda_r^2 + \lambda_\theta^2/r^2}} > 0   \label{eq:pHpsiLOS2}
	\end{array}   \right. 
\end{align}
which means that $\psi^*$ minimizes the cost $-J$; equivalently, it maximizes the final separation $R(t_f)$.

The optimal heading of the Target can be found in a similar way. Let us evaluate 
\begin{align}
  \frac{\partial H}{\partial \phi}=\alpha(1-\lambda_R)\sin\phi - \frac{\alpha}{R}\lambda_\theta \cos\phi = 0. \label{eq:pHphiLOS}
\end{align}
We use the trigonometric identity $\cos^2\phi=1-\sin^2 \phi$ to write \eqref{eq:pHphiLOS} as
\begin{align}
  	  \sin^2\phi^*&=\frac{\lambda_\theta^2}{R^2\big((1\!-\!\lambda_R)^2\! +\! \lambda_\theta^2/R^2\big)}   
\nonumber
\end{align}
and we obtain \eqref{eq:DGsphi}. The expression \eqref{eq:DGcphi} is found in a similar way by setting $\sin^2\phi=1-\cos^2 \phi$ in \eqref{eq:pHphiLOS}. Similarly, we compute  
\begin{align}
\left.
	\begin{array}{l l}
  \frac{\partial^2 H}{\partial \phi^2}&=\alpha(1-\lambda_R) \cos\phi + \frac{\alpha}{R}\lambda_\theta \sin\phi  \\
	&=\frac{\alpha(1-\lambda_R)^2}{\sqrt{(1-\lambda_R)^2 + \lambda_\theta^2/R^2}} + \frac{\alpha\lambda_\theta^2}{R^2\sqrt{(1-\lambda_R)^2 + \lambda_\theta^2/R^2}} > 0   \label{eq:pHphiLOS2}
		\end{array}   \right. 
\end{align}
which means that $\phi^*$ minimizes the cost $-J$; equivalently, it maximizes the final separation $R(t_f)$

The optimal heading $\chi^*$ is characterized in a similar way. We differentiate the Hamiltonian \eqref{eq:HamiltonianDG} in $\chi$ and set the derivative to $0$
\begin{align}
\left.
	\begin{array}{l l}
  \frac{\partial H}{\partial \chi}&=(1-\lambda_R)\sin(\theta-\chi)+\lambda_r\sin\chi  \\
 &~~-\frac{\lambda_\theta}{R}\cos(\theta-\chi)+\frac{\lambda_\theta}{r}\cos\chi = 0.   \label{eq:DGpHchi}
		\end{array}   \right. 
\end{align}
Using the trigonometric identities:
\begin{align}
	\sin(\theta-\chi)=\sin\theta\cos\chi-\cos\theta\sin\chi   \label{eq:sinsum}   \\
	\cos(\theta-\chi)=\cos\theta\cos\chi+\sin\theta\sin\chi       \label{eq:cossum}
\end{align}
we can write \eqref{eq:DGpHchi} as follows:
\begin{align}
 \left.
	\begin{array}{l l}
 \big((1-\lambda_R)\sin\theta-\frac{\lambda_\theta}{R}\cos\theta+\frac{\lambda\theta}{r}\big)\cos\chi  \\
=\big((1-\lambda_R)\cos\theta+\frac{\lambda_\theta}{R}\sin\theta-\lambda_r\big)\sin\chi.  \label{eq:DGcossin}
		\end{array}   \right. 
\end{align}
We now use the trigonometric identity $\cos^2\chi=1-\sin^2\chi$ to obtain 
\begin{align}
   \sin^2\chi^*=  \frac{\chi_s^2}
		{\chi_s^2+\chi_c^2}  \nonumber
\end{align}
and \eqref{eq:DGschi} follows. The expression \eqref{eq:DGcchi} is found in a similar way by setting $\sin^2\chi=1-\cos^2 \chi$ in \eqref{eq:DGcossin}.

In order to guarantee that the Attacker optimal control maximizes the objective $-J$ we evaluate the second partial derivative of the Hamiltonian with respect to the Attacker control input.
\begin{align}
 \left.
	\begin{array}{l l}
  \frac{\partial^2 H}{\partial \chi^2}&=-(1-\lambda_R)\cos(\theta-\chi)+\lambda_r\cos\chi  \\
	&~~-\frac{\lambda_\theta}{R}\sin(\theta-\chi)-\frac{\lambda_\theta}{r}\sin\chi.   \label{eq:DGp2Hchi}
			\end{array}   \right. 
\end{align}
Inserting the expressions \eqref{eq:sinsum} and \eqref{eq:cossum} into eq. \eqref{eq:DGp2Hchi} we obtain the following
\begin{align}
  \left.
	\begin{array}{l l}
	\frac{\partial^2 H}{\partial \chi^2}&= -\big((1-\lambda_R)\sin\theta-\frac{\lambda_\theta}{R}\cos\theta+\frac{\lambda\theta}{r}\big)\sin\chi \\
	&~~- \big((1-\lambda_R)\cos\theta+\frac{\lambda_\theta}{R}\sin\theta-\lambda_r\big)\cos\chi <0.
		\end{array}   \right. 
\end{align}
Therefore, the solutions \eqref{eq:DGschi} and \eqref{eq:DGcchi} maximize the objective $-J$, which is equivalent to minimize the terminal separation $R(t_f)$. $\square$

The expressions for the optimal heading angles \eqref{eq:DGspsi}-\eqref{eq:DGcchi} are used to numerically solve the Two-Point Boundary Value Problem (TPBVP) \eqref{eq:dynR_DG}-\eqref{eq:costDG}, \eqref{eq:costateRDG}-\eqref{eq:DGcchi}. The numerical solution is found by substituting the optimal control headings into the state equations \eqref{eq:dynR_DG}-\eqref{eq:dynD_DG}, and the co-state equations \eqref{eq:costateRDG}-\eqref{eq:costateDG}, with the terminal conditions given by \eqref{eq:DGtc}. 

\section{Critical Speed Ratio for Target Survival}  \label{sec:Minspeed}
\begin{figure}
	\begin{center}
		\includegraphics[width=8.4cm,height=6cm,trim=.5cm .4cm 1.4cm .4cm]{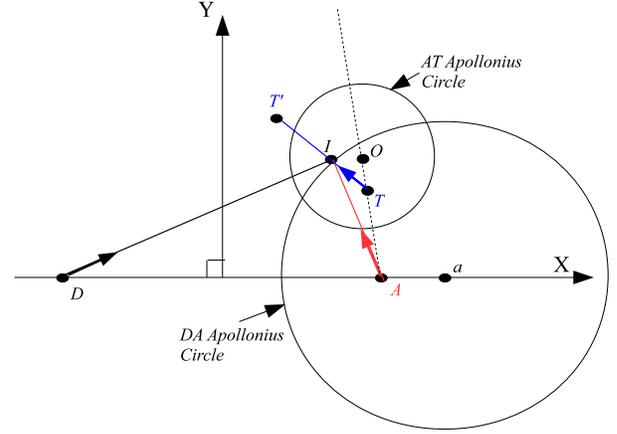}
	\caption{Target-Attacker-Defender scenario for $\gamma<1$}
	\label{fig:Scenario}
	\end{center}
\end{figure}

In this section we consider point capture, that is, the separation $r$ has to satisfy $r(t_f)\rightarrow 0$ in order for the Defender to capture the Attacker. The scenario is illustrated in Fig. \ref{fig:Scenario}.
We consider the rotating reference frame anchored on the Attacker and the Defender.  
In Fig. \ref{fig:Scenario} the points $A$ and $D$ represent the positions of the Attacker and the Defender, respectively. A Cartesian frame is attached to the points $A$ and $D$ in such a way that the extension to infinity of $\overline{AD}$ in both directions represents the $X$-axis and the orthogonal bisector of $\overline{AD}$ represents the $Y$-axis. The positions of the three agents in this frame are $T=(x_T,y_T)$, $A=(x_A,0)$, and $D=(-x_A,0)$.

With respect to Fig. \ref{fig:Scenario} we note that the Attacker aims at minimizing the distance between the Target at the time instant when the Defender intercepts the Attacker, point $T'$, and point $I$, where the Defender intercepts the Attacker. The points $T$ and $T'$ represent the initial and terminal positions of the Target, respectively. 

Define $\gamma=1/\beta=V_A/V_D$. When $\gamma<1$ the Defender will intercept the Attacker at some point $I=(x_I,y_I)$ that lies on the Apollonius circle defined by the Defender and the Attacker separation and the speed ratio $\gamma$. The center of the $DA$-based Apollonius circle is at ($a,0$), where
\begin{align}
	  a=\frac{1+\gamma^2}{1-\gamma^2}x_A \label{eq:a}
\end{align}
and the radius of the $DA$ Apollonius circle is
\begin{align}
	  r_A=\frac{2\gamma}{1-\gamma^2}x_A. \label{eq:rA}
\end{align}
When the Target is inside the $DA$ Apollonius circle, its speed needs to be high enough in order to exit from the $DA$ Apollonius circle before being captured by the Attacker. If the Target is able to exit the $DA$ Apollonius circle then the Defender will be able to assist the Target to escape, by intercepting the Attacker who is on route to the Target. 

\begin{figure}
	\begin{center}
		\includegraphics[width=8.4cm,height=6.5cm,trim=.3cm .1cm .5cm .1cm]{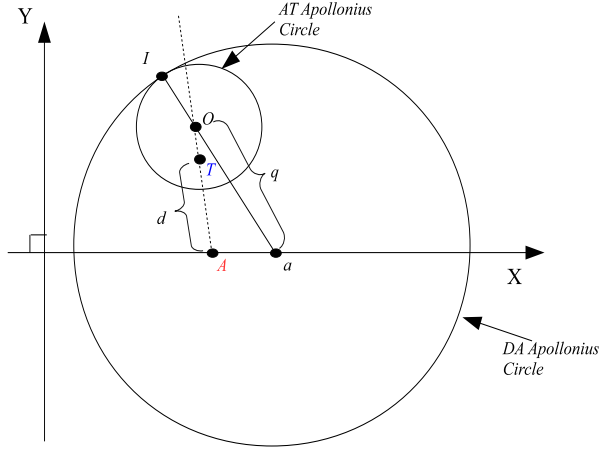}
	\caption{Determination of $\bar{\alpha}$}
	\label{fig:Tangent}
	\end{center}
\end{figure}

\begin{proposition}
Given the speed ratio $\gamma=V_A/V_D<1$, the critical speed ratio $\bar{\alpha}$ is a function of the positions of the Target and the Attacker and is given by
\begin{align}
     \bar{\alpha}=\frac{\gamma\sqrt{(x_A+x_T)^2+y_T^2}-\sqrt{(x_A-x_T)^2+y_T^2}}{2\gamma x_A}.    \label{eq:alphasol}  
\end{align}
\end{proposition}
\textit{Proof}. 
In order to determine the minimum speed ratio, $\bar{\alpha}$, that guarantees Target survival we consider a second Apollonius circle defined by the Attacker and the Target using the Target/Attacker speed ratio $\alpha$. Thus, a solution to the differential game exists if and only if the $AT$ Apollonius circle, which is based on the segment $\overline{AT}$ and the speed ratio $\alpha$, intersects the $DA$ Apollonius circle, the one based on the the segment $\overline{DA}$ and the speed ratio $\gamma$. The lower limit $\bar{\alpha}$ on the speed ratio $\alpha$, that is, $\bar{\alpha}<\alpha<1$, corresponds to the case where the $AT$ Apollonius circle is tangent to the $DA$ Apollonius circle, see Fig. \ref{fig:Tangent}. Note that if the speed ratio $\alpha\geq 1$ the Target always escapes and there is no need for a Defender missile, that is, there is no target defense differential game in the first place.

The Attacker's initial position, the Target's initial position, and the center $O$ of the $AT$ Apollonius circle are collinear and lie on the dotted line in Fig. \ref{fig:Tangent} which can be represented as
\begin{align}
	  y=-\frac{y_T}{x_A-x_T}x + \frac{x_Ay_T}{x_A-x_T}. \nonumber
\end{align}
The geometry of the second Apollonius circle is as follows: The center of the circle, denoted by $O$, is at a distance of $\frac{\alpha^2}{1-\alpha^2}d$ from $T$ and its radius is $r_O=\frac{\alpha}{1-\alpha^2}d$, where $d$ is the distance between $A$ and $T$ and is given by
\begin{align}
   d=\sqrt{(x_A-x_T)^2+y_T^2}. 
\end{align}
Hence, the following holds
\begin{align}
	  & \Big(\frac{x_Ty_T}{x_A-x_T}-\frac{y_T}{x_A-x_T}x_0\Big)^2 + (x_0-x_T)^2  \nonumber  \\
		&=\frac{\alpha^4}{(1-\alpha^2)^2}[(x_A-x_T)^2+y_T^2]   \nonumber  
\end{align}
and we calculate the coordinates of the center of the second Apollonius circle
\begin{align}
   \left.
	 \begin{array}{l l}
     x_O=\frac{1}{1-\alpha^2}x_T-\frac{\alpha^2}{1-\alpha^2}x_A \\
		y_O=\frac{1}{1-\alpha^2}y_T.
	\end{array}  \label{eq:circCenter}  \right.
\end{align}
From Fig. \ref{fig:Tangent} we can see that the three points $a$, $O$, and $I$ are collinear, where $I$ represents the tangent point where both circles meet. Thus, we have the following relationship
\begin{align}
	  r_O=r_A-q  \label{eq:rOrA}
\end{align}
where $q=\sqrt{(a-x_O)^2+y_O^2}$. Eq. \eqref{eq:rOrA} can be written as follows
\begin{align}
	  (a-x_O)^2+y_O^2=(r_A-\frac{\alpha}{1-\alpha^2}d)^2.  \label{eq:alpb1}
\end{align}
Eq. \eqref{eq:alpb1} can be expressed in terms of the known positions ($x_T$, $y_T$, and $x_A$), the known speed ratio $\gamma$, and the variable we aim to solve for, which is $\alpha$. After a few steps we obtain the following quartic equation in $\alpha$ 
\begin{align}
	 \left.
	 \begin{array}{l l}
  \frac{4\gamma^2}{1-\gamma^2}x_A^2\alpha^4 + \frac{4\gamma d}{1-\gamma^2}x_A\alpha^3 \\ 
	- \Big(\frac{4\gamma^2}{1-\gamma^2}x_A^2+\frac{4\gamma^2}{1-\gamma^2}x_Ax_T-d^2\Big)\alpha^2 \\
	- \frac{4\gamma d}{1-\gamma^2}x_A\alpha +\frac{4\gamma^2}{1-\gamma^2}x_Ax_T -d^2 =0  
	\end{array}  \right.
\end{align}
which can be factored out to obtain the quadratic polynomials
\begin{align}
 \left.
	 \begin{array}{r r}
	\big(4\gamma x_A\alpha(\gamma x_A\alpha +d) + (1-\gamma^2)d^2 - 4\gamma^2 x_Ax_T\big) \big(\alpha^2-1\big) \\  
	= 0.  
		\end{array}    \right.  \nonumber
\end{align}
The solutions $\alpha=\pm 1$ are irrelevant to the differential game under analysis. Thus, the critical speed ratio $\bar{\alpha}$ is given by the positive solution of the quadratic equation
\begin{align}
	4\gamma^2 x_A^2\alpha^2 + 4\gamma x_Ad\alpha + (1-\gamma^2)d^2 - 4\gamma^2 x_Ax_T  = 0   \nonumber
\end{align}
which is given by \eqref{eq:alphasol}.              $\square$

In the particular case where $y_T=0$, the critical speed ratio is given by
\begin{align}
     \bar{\alpha}=\frac{x_A(\gamma-1)+x_T(\gamma+1)}{2\gamma x_A}.   \label{eq:alpb2}
\end{align}
Further, the initial Target position ($x_T,0$) for which the critical $\bar{\alpha}$ is equal to zero can be obtained from \eqref{eq:alpb2}
\begin{align}
     0=\frac{x_A(\gamma-1)+x_T(\gamma+1)}{2\gamma x_A}   \nonumber  \\
		\Rightarrow x_T=\frac{(1-\gamma)x_A}{1+\gamma}
\end{align}
which is equivalent to 
\begin{align}
		a-r_A=\frac{1+\gamma^2}{1-\gamma^2}x_A-\frac{2\gamma}{1-\gamma^2}x_A = \frac{(1-\gamma)x_A}{1+\gamma},  \nonumber
\end{align}
as expected.

\section{Optimal Strategies} \label{sec:analysis}
When $\gamma<1$ the Defender will intercept the Attacker at some point $I=(x_I,y_I)$ that lies on the $DA$ Apollonius circle. 
Notice that all points outside the $DA$ Apollonius circle can be reached by the Defender before the Attacker does; similarly, all points inside the same circle can be reached by the Attacker before the Defender does.

\subsection{Target Starts Outside of $DA$ Apollonius Circle}  
In the case where the Target is initially outside the $DA$ Apollonius circle it can be clearly seen that the Defender can help the Target regardless of the speed ratio $0<\alpha<1$ because the Attacker cannot reach the Target before the Defender does. In other words, the critical speed ratio in this case is $\bar{\alpha}=0$.
In the case where the Target is outside the $DA$ Apollonius circle, the Target chooses point $v$ on the $DA$ Apollonius circle (in order to run away from that point) and the Attacker chooses his aimpoint $u$ on the same circle. Additionally, the Defender tries to intercept the Attacker by choosing his aimpoint $w$, also on the $DA$ Apollonius circle. The Target, the Defender, and the Attacker are faced with the minmax optimization problem: $\min_u \max_{v,w} J(u,v,w)$, where $J(u,v,w)=S$ and $S$ represents the distance between the Target terminal position $T'$ and the point on the $DA$ Apollonius circle where the Attacker is intercepted by the Defender. 

The Defender helps the Target to escape by intercepting the Attacker at the point $u$ on the $DA$ Apollonius circle. Therefore, the Defender's optimal policy is $w^*(u,v)=u$ in order to guarantee interception of the Attacker.
Since the Defender's optimal policy is $w^*=u$ we have that the decision variables $u$ and $v$ jointly determine $J(u,v)$, where $J(u,v)=S$ and $S$ represents the distance between the Target terminal position $T'$ and the point $I=u$ on the $DA$ Apollonius circle where the Attacker is intercepted by the Defender. 

\begin{proposition}  \label{prop:outside}
Given the cost/payoff function $J(u,v)$, the solution $u^*$ and $v^*$ of the optimization problem $ \min_{u} \ \max_{v} J(u,v)$ is such that
\begin{align}
    u^*=v^*.  \nonumber
\end{align}
Moreover, when the Target is outside the $DA$ Apollonius circle, the Target's strategy is $v^*(u)=\arg \max_v J(u,v)=u$ so that it suffices to solve the optimization problem 
\begin{align}
  \left.
	 \begin{array}{l l}
		\min_{x_I,y_I} J(x_I,y_I)
		\\
		\text{subject to} \ \ (a-x_I)^2+y_I^2 = r_A^2
		\end{array}  \label{eq:JxyProb}  \right.
\end{align}
where
\begin{align}
  \left.
	 \begin{array}{l l}
		J(x_I,y_I)\!\!\!&=\sqrt{(x_I-x_T)^2+(y_I-y_T)^2} \\
		  &~~+ \alpha\sqrt{(x_A-x_I)^2+y_I^2}.  \\
		\end{array}  \label{eq:Jxy}  \right.
\end{align}
\end{proposition}
\begin{flushright}
$\square$
\end{flushright}

\begin{figure}
	\begin{center}
		\includegraphics[width=8.4cm,height=6.5cm,trim=.1cm .1cm .4cm .1cm]{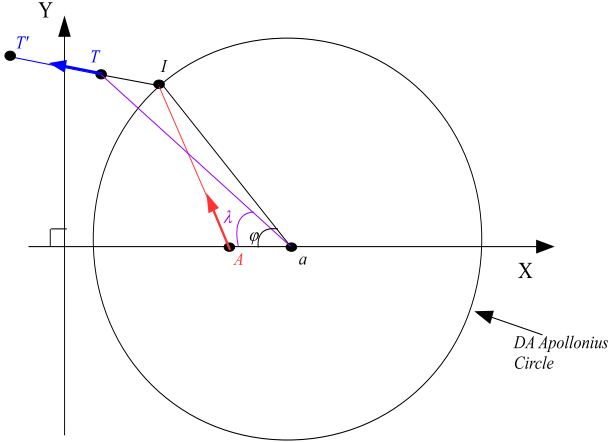}
	\caption{Optimal strategy}
	\label{fig:Optimal}
	\end{center}
\end{figure}

The Attacker chooses the optimal coordinates $(x_I,y_I)$ of point $I$ that minimize the final separation $J(x_I,y_I)=\overline{IT}+\overline{TT'}$ see Fig. \ref{fig:Optimal}.

One way to formulate this problem is as shown in Proposition \ref{prop:outside}. The equality constraint can be used in the cost to write $J(x_I)$, that is, to write the cost in terms of only one variable. For instance, the first derivative $\frac{\partial J(x_I)}{\partial x_I}=0$ results in a sixth order equation in $x_I$. 

\begin{theorem}
The optimal interception point $I$ that minimizes \eqref{eq:Jxy} has polar coordinates $I=(\varphi^*,r_A)$ with respect to the center of the $DA$ Apollonius circle denoted by $a$, where $\varphi^*$ is the solution of the sixth order complex exponential equation
\begin{align}
  \left.
	 \begin{array}{l l}
  \frac{Nr_A}{l}\big(1-\frac{N}{\alpha^2Ml}\big)e^{6i\varphi} \\
	+ \big((\frac{N}{\alpha Ml})^2(r_A^2+M^2)-r_A^2-N^2\big)e^{5i\varphi} \\
	+ Nr_A\big(\frac{N}{\alpha^2Ml^2}(2l^2-1)+l-\frac{2}{l}  \big) e^{4i\varphi}  \\
	+2\big(r_A^2+N^2-(\frac{N}{\alpha M})^2(r_A^2+M^2)\big)e^{3i\varphi} \\
	+Nr_A\big(\frac{N}{\alpha^2M}(2-l^2)-2l+\frac{1}{l}\big) e^{2i\varphi} \\
	+ \big((\frac{Nl}{\alpha M})^2(r_A^2+M^2)-r_A^2-N^2\big)e^{i\varphi} \\
	+ Nr_Al\big(1-\frac{Nl}{\alpha^2M}\big) = 0
	\end{array}  \label{eq:fderJ2}  \right.
\end{align}
that minimizes the cost 
\begin{align}
 \left.
	 \begin{array}{l l}
  J(\varphi)\!\!\!&=\sqrt{r_A^2 + N^2 - 2Nr_A\cos(\varphi -\lambda)} \\
	&~~ +\alpha\sqrt{r_A^2 + M^2 - 2Mr_A\cos\varphi}  
\end{array}   \right.  \label{eq:Jvp}
\end{align}
where $l=e^{i\lambda}$, $M=\frac{2\gamma^2}{1-\gamma^2}x_A$ represents the distance between the points $A$ and $a$, and $N=\sqrt{(a-x_T)^2+y_T^2}$ represents the distance between the points $a$ and $T$.
\end{theorem}

\textit{Proof}. An alternative and more compact algebraic equation to directly solving \eqref{eq:Jxy} can be obtained by searching for the optimal angle $\varphi$ that minimizes the same cost see Fig. \ref{fig:Optimal}.

It can be seen that, by varying the angle $\varphi$, the point $I$ moves along the circumference of the $DA$ Apollonius circle. In order to write an equivalent expression to \eqref{eq:Jxy}, but only in terms of $\varphi$, we consider the two triangles $\Delta aAI$ and $\Delta aTI$.

The distance $\overline{TT'}$ is proportional to the distance  $\overline{AI}$. The distance $\overline{AI}$ changes as the angle $\varphi$ takes different values. However, the distance $\overline{aI}=r_A$ and the distance $\overline{aA}=\frac{2\gamma^2}{1-\gamma^2}x_A$ are fixed. Similarly, the distance $\overline{IT}$ changes in terms of the angle $\varphi$, but the distance $\overline{aI}$, the distance $\overline{aT}$, and the angle $\lambda$ are fixed.
Then, the cost \eqref{eq:Jxy} can be written as in \eqref{eq:Jvp}

The first derivative of \eqref{eq:Jvp} is
\begin{align}
\left.
	 \begin{array}{l l}
  \frac{d J(\varphi)}{d\varphi}\!\!\! &= \frac{N\sin(\varphi -\lambda)}{\sqrt{r_A^2 + N^2 - 2Nr_A\cos(\varphi -\lambda)}} \\
	&~~+ \frac{\alpha M \sin\varphi}{\sqrt{r_A^2 + M^2 - 2Mr_A\cos\varphi}}.  \label{eq:fdervp}
\end{array}   \right.
\end{align}
Setting \eqref{eq:fdervp} equal to zero we obtain
\begin{align}
\left.
	 \begin{array}{l l}
 \frac{N^2\sin^2(\varphi -\lambda)}{r_A^2 + N^2 - 2Nr_A\cos(\varphi -\lambda)} =
	 \frac{\alpha^2 M^2 \sin^2\varphi}{r_A^2 + M^2 - 2Mr_A\cos\varphi}.  \label{eq:fdervp2}
\end{array}   \right.
\end{align}
In order to solve for the angle $\varphi$ we 
use the complex exponential $e^{i\varphi}$ to obtain
\begin{align}
  \left.
	 \begin{array}{l l}
  \frac{N^2}{4}\big(e^{i(\varphi-\lambda)}\!-\! e^{-i(\varphi-\lambda)}\big)^2 \big(r_A^2\!+\!M^2\!-\! Mr_A(e^{i\varphi}\!+\!e^{-i\varphi})\big) = \\
	\frac{\alpha^2M^2}{4}\big(e^{i\varphi}\!-\!e^{-i\varphi}\big)^2 \big(r_A^2\!+\!N^2\!-\!Nr_A(e^{i(\varphi-\lambda)}+e^{-i(\varphi-\lambda)})\big).
	\end{array}  \label{eq:fderJ}  \right.
\end{align}
After some manipulation we obtain a sixth order polynomial equation in $e^{i\varphi}$ as it is shown in \eqref{eq:fderJ2}.

The polynomial in \eqref{eq:fderJ2} has complex coefficients. The six solutions of \eqref{eq:fderJ2} are complex, in general, of the form $e^{i\varphi}=\cos\varphi+i\sin\varphi$. Thus, the angle $\varphi$ can be directly obtained. In the worst case, we only need to test the six angles in the cost function to determine the optimal solution $\varphi^*$.  $\square$

\subsection{Target Starts Inside of $DA$ Apollonius Circle} 
In the case where the Target is inside the $DA$ Apollonius circle, the Target chooses his aimpoint $v$ on the $DA$ Apollonius circle and the Attacker chooses his aimpoint $u$ on the same circle. Additionally, the Defender tries to intercept the Attacker by choosing his aimpoint $w$, also on the $DA$ Apollonius circle. The Target, the Defender, and the Attacker are faced with the maxmin optimization problem: $\max_{v,w} \min_u J(u,v,w)$, where $J(u,v,w)=S$ and $S$ represents the distance between the Target terminal position $T'$ and the point on the $DA$ Apollonius circle where the Attacker is intercepted by the Defender. 

\begin{figure}
	\begin{center}
		\includegraphics[width=8.4cm,height=6.5cm,trim=.4cm .2cm 1.4cm .2cm]{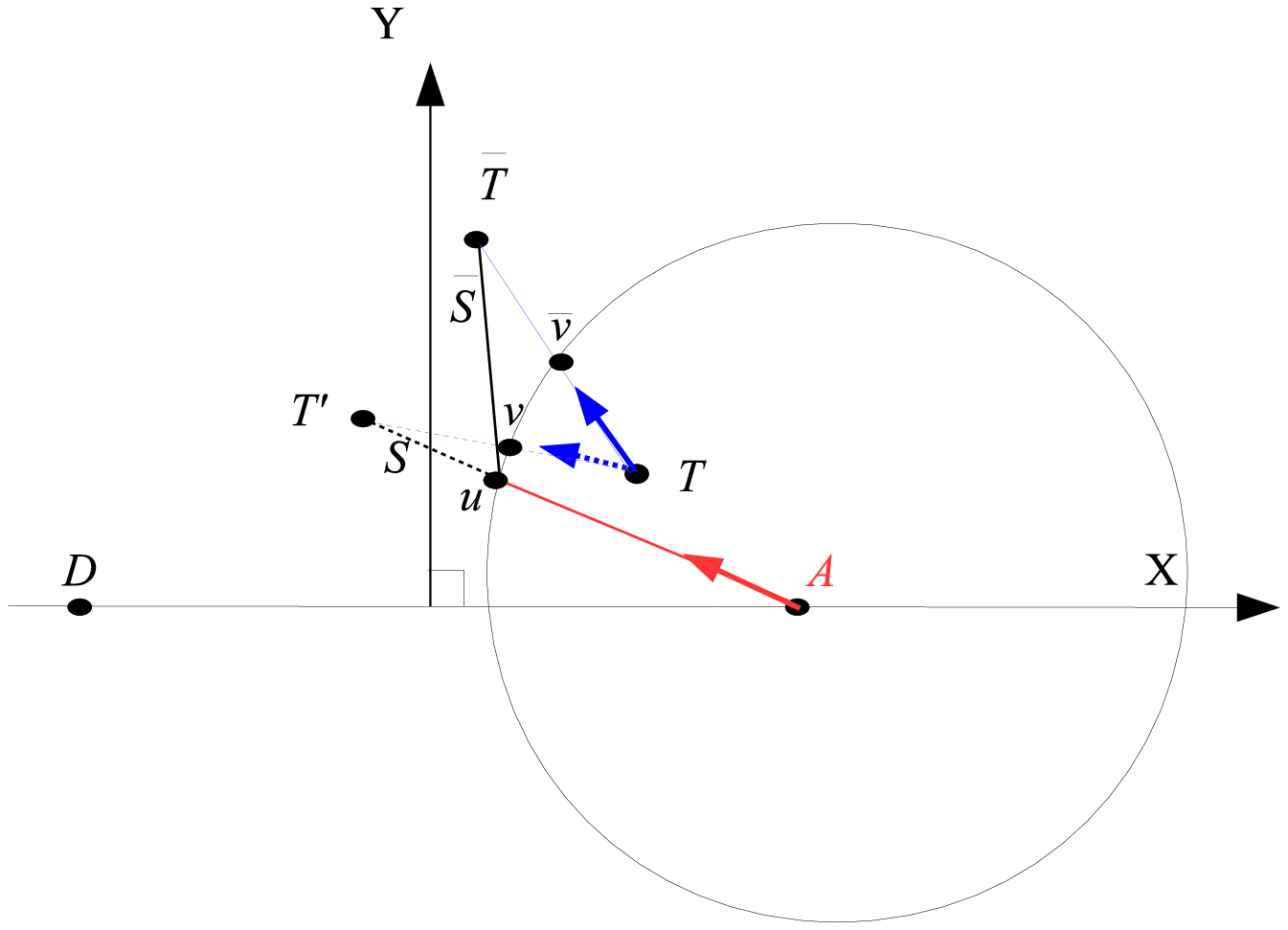}
	\caption{maxmin optimization problem}
	\label{fig:maxmin}
	\end{center}
\end{figure}

The Defender's optimal policy is $w^*(u,v)=u$ in order to guarantee interception of the Attacker and the decision variables $u$ and $v$ jointly determine $J(u,v)$.
Now, let us analyze the possible strategies. If the Target chooses $v$, the Attacker will respond and choose $u$. If $u\neq v$ the Target would correct his decision and choose some $\bar{v}$ such that $\bar{S}>S$ as shown in Fig. \ref{fig:maxmin}. In general, choosing $u\neq v$ is detrimental to the Attacker since the resulting cost will increase. Thus, it is clear that the Attacker should aim at the point $v$ which is chosen by the Target.

\begin{proposition}
Given the cost/payoff function $J(u,v)$, the solution $u^*$ and $v^*$ of the optimization problem $ \max_{v} \ \min_{u} J(u,v)$ is such that
\begin{align}
    u^*=v^*.   \nonumber
\end{align}
Moreover, when the Target is inside the $DA$ Apollonius circle, the Attacker's strategy is $u^*(v)=\arg \min_u J(u,v)=v$ so that it suffices to solve the optimization problem 
\begin{align}
  \left.
	 \begin{array}{l l}
		\max_{x_I,y_I} J(x_I,y_I)
		\\
		\text{subject to} \ \ (a-x_I)^2+y_I^2 = r_A^2
		\end{array}  \label{eq:JxyProb2}  \right.
\end{align}
where
\begin{align}
  \left.
	 \begin{array}{l l}
		J(x_I,y_I)\!\!\!&=\alpha\sqrt{(x_A-x_I)^2+y_I^2} \\
		  &~~- \sqrt{(x_I-x_T)^2+(y_I-y_T)^2}.  \\
		\end{array}  \label{eq:Jxy2}  \right.
\end{align}
\end{proposition}
\begin{flushright}
$\square$
\end{flushright}

The difference between the Target being inside or outside the $DA$ Apollonius circle is not only the sign in the cost function but the Target and Attacker strategies. In the case treated in this subsection, the Target chooses the coordinates $(x_I,y_I)$ that maximize the final separation $J(x_I,y_I)$ and the Attacker follows the Target's decision. Additionally, when the Target is inside the $DA$ Apollonius circle its critical speed is $\bar{\alpha}>0$ (given by \eqref{eq:alphasol}); when the Target is outside the same circle then its critical speed is $\alpha=0$.

\begin{theorem}
The optimal interception point $I$ that maximizes \eqref{eq:Jxy2} has polar coordinates $I=(\varphi^*,r_A)$ with respect to the center of the $DA$ Apollonius circle denoted by $a$, where $\varphi^*$ is the solution of the sixth order complex exponential equation \eqref{eq:fderJ2}
that maximizes the cost 
\begin{align}
 \left.
	 \begin{array}{l l}
  J(\varphi)&= \alpha\sqrt{r_A^2 + M^2 - 2Mr_A\cos\varphi}  \\ 
	&~~-\sqrt{r_A^2 + N^2 - 2Nr_A\cos(\varphi -\lambda)}   
\end{array}   \right.  \label{eq:Jvpin}
\end{align}
where $l=e^{i\lambda}$, $M$ represents the distance between the points $A$ and $a$, and $N$ represents the distance between the points $a$ and $T$.
\end{theorem}

\textit{Proof}.
The cost \eqref{eq:Jxy2} can be written in terms of the angle $\varphi$ as in \eqref{eq:Jvpin}. 
The first derivative of \eqref{eq:Jvpin} is
\begin{align}
\left.
	 \begin{array}{l l}
  \frac{d J(\varphi)}{d\varphi} &= \frac{\alpha M \sin\varphi}{\sqrt{r_A^2 + M^2 - 2Mr_A\cos\varphi}}  \\
	&~~-  \frac{N\sin(\varphi -\lambda)}{\sqrt{r_A^2 + N^2 - 2Nr_A\cos(\varphi -\lambda)}}.  \label{eq:fdervpin}
\end{array}   \right.
\end{align}
Setting \eqref{eq:fdervpin} equal to zero we obtain \eqref{eq:fdervp2} and, consequently, the optimal angle $\varphi^*$ is the solution of \eqref{eq:fderJ2} that maximizes \eqref{eq:Jvpin}. $\square$

\section{Examples}    \label{sec:examples}
\begin{figure}
	\begin{center}
		\includegraphics[width=8.4cm,height=6.5cm,trim=.9cm .6cm .9cm .6cm]{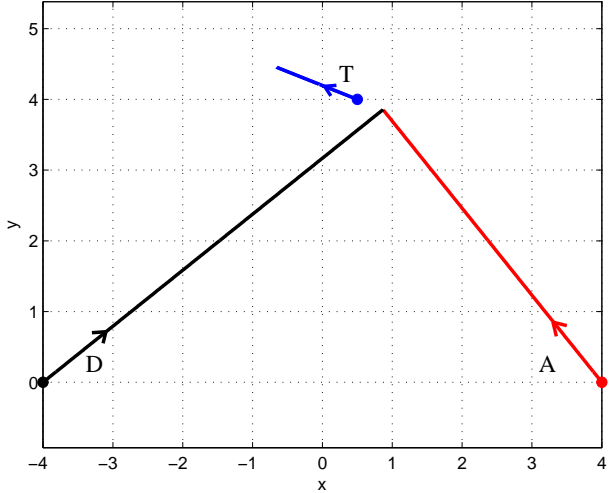}
	\caption{Optimal trajectories in Example 1}
	\label{fig:Ex1}
	\end{center}
\end{figure}

\textit{Example 1. Target is outside the DA Apollonius circle}. The speed ratios are $\alpha=0.25$ and $\gamma=0.8$. The initial conditions of the three agents are given by: $A=(4,0)$, $D=(-4,0)$, and $T=(0.5,4)$. We calculate $a=18.22$ and $r_A=17.78$. The six solutions of \eqref{eq:fderJ2} are given by
\begin{align}
\left.
	 \begin{array}{l l}
  \varphi_1=-2.9596  \\
	\varphi_2=-2.8573  \\
	\varphi_3=0.0001  \\
	\varphi_4=0.0001  \\
	\varphi_5=0.2254  \\
	\varphi_6=0.2186 . 
	\end{array}  \nonumber  \right.
\end{align}
By evaluating these solutions using \eqref{eq:Jvp} we have that the optimal solution is $\varphi^*=0.2186$, which yields $I^*=(0.8676,3.8555)$. The trajectories are shown in Fig. \ref{fig:Ex1}. Note that the same trajectories and optimal interception point are obtained by using the numerical method from Sec. \ref{sec:Num}.

\begin{figure}
	\begin{center}
		\includegraphics[width=8.4cm,height=6.5cm,trim=.9cm .6cm .9cm .6cm]{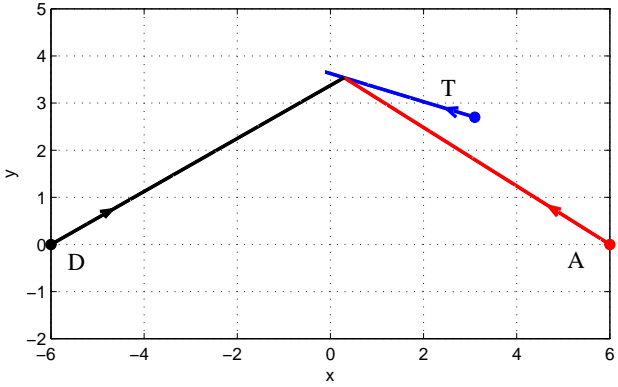}
	\caption{Optimal trajectories in Example 2}
	\label{fig:Ex2}
	\end{center}
\end{figure}

\textit{Example 2. Target is inside the DA Apollonius circle}. The speed ratios are $\alpha=0.5$ and $\gamma=0.93$. The initial conditions of the three agents are given by: $A=(6,0)$, $D=(-6,0)$, and $T=(3.1,2.7)$. In this case we calculate $a=82.823$ and $r_A=82.605$. Since the Target is initially inside the $DA$ Apollonius circle, the critical speed ratio $\bar{\alpha}$ is greater than zero. We can use eq. \eqref{eq:alphasol} to find the exact value of the critical speed ratio which is $\bar{\alpha}=0.436$. Thus, the value $\alpha=0.5>\bar{\alpha}$ guarantees the Target's escape. Now we can search for the optimal angle $\varphi^*$ that solves the differential game.
The six solutions of eq. \eqref{eq:fderJ2} are given by
\begin{align}
\left.
	 \begin{array}{l l}
  \varphi_1=-3.0752  \\
	\varphi_2=-3.1189  \\
	\varphi_3=-0.0014  \\
	\varphi_4=-0.0014  \\
	\varphi_5=0.0277  \\
	\varphi_6=0.0429 . 
	\end{array}  \nonumber  \right.
\end{align}
By evaluating these solutions using eq. \eqref{eq:Jvpin} we have that the optimal solution is $\varphi^*=0.0429$, which yields $I^*=(0.293,3.539)$. The trajectories are shown in Fig. \ref{fig:Ex2}. Note that the same trajectories and optimal interception point are obtained by using the numerical method from Sec. \ref{sec:Num}.

\textit{Example 3. Robustness to unknown Attacker guidance law}. A very important characteristic of the cooperative guidance laws for the active target defense differential game as discussed in this paper is that the solution given by the sixth order equation \eqref{eq:fderJ2} is a closed-loop interception strategy that is robust to unknown Attacker guidance laws. This means that if the Attacker does not follow its optimal policy and uses a different guidance law that is unknown to the Target-Defender team then the Target and the Defender (having current measurements of the Attacker' position) are able to solve \eqref{eq:fderJ2} and continuously update their cooperative interception strategy, thus increasing the $T-A$ separation at interception time. 

\begin{figure}
	\begin{center}
		\includegraphics[width=8.4cm,height=6.5cm,trim=.9cm .2cm .9cm .2cm]{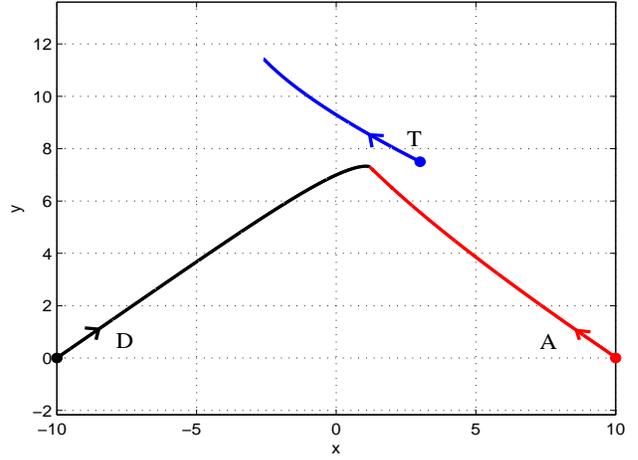}
	\caption{Trajectories in Example 3}
	\label{fig:Ex4}
	\end{center}
\end{figure}

Let the initial positions of the three agents be: $A=(10,0)$, $D=(-10,0)$, and $T=(3,7.5)$. The speed ratios are $\alpha=0.6$ and $\gamma=0.85$. The Attacker implements $PN$ guidance law with navigation constant N=3. However, this information is unknown to the Target-Defender team and they are only able to measure the current position of the Attacker, $A=(x_A(t),y_A(t))$. By continuously updating their headings, the Target-Defender team are able to defeat the Attacker, that is, the Defender intercepts the Attacker and the Target escapes being captured by the Attacker. The trajectories for this example are shown in Fig. \ref{fig:Ex4}. The final separation between Target and Attacker is $R(t_f)=5.609>J^*$. As expected, the final separation is more than if the Attacker played optimally. When the Attacker plays optimally the cost/payoff is $J^*=5.373$.

\section{Conclusions} \label{sec:concl}
A numerical and an analytical solution to the active target defense differential game with a fast Defender were presented in this paper. The numerical solution is based on the Pontryagin's Maximum Principle applied to the differential game and a TPBVP is solved numerically. The analytical approach hinges on the solutions of a sixth-order polynomial equation that provides the optimal interception point's coordinates, hence it provides the optimal headings for the players. This result comes with an expected increase in complexity compared to the case where both missiles are restricted to have the same speed \cite{Pachter14Allerton}. In that case, the solution of the differential game required the rooting of a fourth-order polynomial.


\end{document}